\documentclass{article}
\usepackage[utf8]{inputenc}
\usepackage[margin = 1.2 in]{geometry}

\usepackage{graphicx}
\usepackage{amsfonts}
\usepackage{amsmath}
\usepackage{caption}
\usepackage{subcaption}
\usepackage{cite}
\usepackage{mathtools}
\usepackage{hyperref}

\usepackage{authblk}

\hypersetup{colorlinks,linkcolor={blue},citecolor={red},urlcolor={red}}

\title{\textbf{Reduced-Order Model of the Russian Service Module via Loewner Framework}}
\author[1]{\textbf{Sanwar Alam}}
\author[2,3]{\textbf{Mohammad N. Murshed}}
\affil[1]{Department of Electrical and Computer Engineering, North South University}
\affil[2]{Department of Mathematics and Physics, North South University}
\affil[3]{Department of Mathematics, University of Tennessee at Chattanooga}

\begin{document}

\maketitle

\begin{abstract}
    \textbf{Loewner framework is a technique that uses frequency response data to construct a reduced order model of a given system. In the past, it has been employed in many different synthetic problems and applications like beams. In this work, we exploit the tool on data pertaining to the structural vibrations in the Russian Service module. As per our analysis, the Loewner model performs as good as the original system.}
\end{abstract}

\vspace{0.5 cm}

\section{Introduction}
Model Reduction (MOR) has been an interesting aspect of computational science. The idea of MOR is to represent a large system (that has many equations in its original form) in just a few equations while retaining the structure of the dynamics. There are quite a few data-driven model reduction tools like proper orthogonal decomposition \cite{chatterjee2000introduction}, balanced proper orthogonal decomposition \cite{rowley2005model}, Eigensystem Realization Algorithm \cite{murshed2021towards}, dynamic mode decomposition \cite{murshed2019time}, and many more. Our study is centered on a recently developed model reduction technique known as the Loewner framework \cite{antoulas2017tutorial} that works solely based on the frequency response data of the system.\\ \\
The Loewner framework extracts essential underlying system structure by collecting the necessary information from large amount of data. This technique facilitates the derivation of state-space models directly from measured input/output data \cite{antoulas2016loewner}. %Such framework was first developed for data-driven reduction of linear system \cite{antoulas2017tutorial}. 
Reduction of linear switched system was recently generalised by extending this framework allowing different dimensions of the reduced state-space in different modes \cite{gosea2018data}. Relying on Loewner framework, Ionita and Antoulas illustrated a new model of parametrized systems \cite{ionita2014data}, which is able to choose separate reduced orders for each parameter based on the rank of the Loewner matrices. They have also applied the tool on frequency response measurements with \(4\) inputs and \(4\) outputs and created a low order model based on the first few dominant singular values of the Loewner matrix. In \cite{gosea2017approximation}, the authors showed improvement in performance (in terms of approximation) exhibited by the Loewner method compared to modal truncation for a specific target range of frequency. Here, they considered the frequency response of the original beam model and used the framework to find a linear model. %After that, they constructed samples and a bunch of number of interpolation points by using a particular interest range of frequency. Also, it was decided to compress from certain number of singular values by seeing the fast decay.
\cite{gosea2016stability} made an extension of the Loewner model by developing reliable and robust post processing methods to adjust the poles so to yield a reduced model to deal with situations where the standard Loewner model happens to lack stability. \\ \\
In this project, we exclusively apply Loewner framework on a particular application namely the vibrations that occur as an airship docks at the international space station. Thorough systems analysis is done and a reduced order model is subsequently created. The rest of the paper is organized as follows. The fundamentals of the Loewner model are discussed in Section \ref{BG} and the numerical results are shown in Section \ref{NR}. A conclusion is drawn in Section \ref{CFW}.
%They illustrated this by two examples. First one was PEEC model with one input and one output. After truncation, they found of the Loewner model has two unstable poles. Then these poles were excluded by flipping the antistable poles approximant after using the optimal and sub-optimal stable approximant. For second experiment they consider real life measurements from an industry application. Here they recorded the number of antistable poles and varied a particular number of reduction order for truncation.

\section{Background}
\subsection{System Dynamics and Loewner Framework}
\label{BG}
%Loewner framework has been utilised to identify the system using frequency domain measurement for rational interpolation/approximation problems. 
We consider the linear, time-invariant system, \\ 
\begin{equation}
\begin{cases}
    \boldsymbol{E} \boldsymbol{\dot x}(t) = \boldsymbol{A}\boldsymbol{x}(t) + \boldsymbol{B}\boldsymbol{u}(t),\\
    \boldsymbol{y}(t) = \boldsymbol{C}\boldsymbol{x}(t),
    \end{cases}
\end{equation}
where \(\boldsymbol{x}(t) \in \mathbb{R}^k\), \(\boldsymbol{u}(t) \in \mathbb{R}^p\), \(\boldsymbol{y}(t) \in \mathbb{R}^q\) refer to the states, the inputs, and the outputs, respectively. Note that the size of \(p\) and \(q\) are much smaller than \(k\) and that \(\boldsymbol{E}\) is invertible.\\ \\
In frequency-domain analysis, let \(\boldsymbol{s}_{i}\) be the frequencies and \(\boldsymbol{H}_{i}\) be the transfer function values that satisfy \(\boldsymbol{H}(s) = \boldsymbol{C}(s\boldsymbol{E}-\boldsymbol{A})^{-1}\boldsymbol{B}\). \\ \\
% Assuming the points is given as \(s_i\) in complex plane and corresponding transfer function measurement as \(\boldsymbol{S}_i \in \mathbb{C}^{q\times p}\) can be defined, 
% \begin{equation}
%     \boldsymbol{H}(s_i) = \boldsymbol{S}_i, \qquad i = 1, \dots, N,
% \end{equation}
% which is satisfied by associated transfer function \(\boldsymbol{H}(s) = \boldsymbol{C}(s\boldsymbol{E}-\boldsymbol{A})^{-1}\boldsymbol{B}\). Inputs and outputs is referred as \(p\) and \(q\), and size of \(p\) and \(q\) assumed to be much smaller than \(N\).\\ \\
From given measurements, data is composed of two separated tangential interpolations as Right interpolation data and Left interpolation data, and this tangential direction is defined as altering columns/rows of the identity matrix \cite{lefteriu2009new}. 
%The original matrix \(\boldsymbol{S}_i\) will be generated by this resulting vector data as column and row vectors. \\ 

\noindent Right interpolation data is constructed as,
\begin{equation}
    \{(\lambda_i, r_i, w_i|\lambda_i\: \in\: \mathbb{C}, r_i\: \in\: \mathbb{C}^{m\times1}, w_i\:\in\:\mathbb{C}^{p\times1},i = \overline{1,\rho})\}.
\end{equation}
Or, in other words, \\

\(\boldsymbol{\Lambda}\) = \texttt{diag}\([\lambda_1,...,\lambda_\rho]\in\mathbb{C}^{\rho\times \rho},\) \(\boldsymbol{R}\) = \([r_1,...,r_k]\in\mathbb{C}^{m\times \rho},\) \(\boldsymbol{W}\) = \([w_1,...,w_\rho]\in\mathbb{C}^{\rho\times  p}\). \\

\noindent And left interpolation data is constructed as,
\begin{equation}
    \{(\mu_j, l_j, v_j|\mu_j\: \in\: \mathbb{C}, l_j\: \in\: \mathbb{C}^{1\times p}, v_j\;\in\:\mathbb{C}^{1\times m},j = \overline{1,v})\}.
\end{equation}
Or, in other words, \\

\(\boldsymbol{M}\) = \texttt{diag}\([\mu_1,...,\mu_v]\in\mathbb{C}^{v\times v},\) \(\boldsymbol{L}\) = \(\begin{bmatrix} l_1 \\ \vdots \\ l_v \end{bmatrix} \in \mathbb{C}^{v \times p},\) \(\boldsymbol{V}\) = \(\begin{bmatrix} v_1 \\ \vdots \\ v_v \end{bmatrix} \in \mathbb{C}^{v \times m}\). \\ 

\noindent The Loewner matrix \(\mathbb{L}\) and shifted Loewner matrix \(\mathbb{L}_\sigma\) are defined by the realization of descriptor form \(\boldsymbol{E,A,B,C}\) of transfer function \(\boldsymbol{H}(s) = \boldsymbol{C}(s\boldsymbol{E}-\boldsymbol{A})^{-1}\boldsymbol{B}\), which satisfies the right tangential constraints:
\begin{equation}
    H(\lambda_i)r_i = [\boldsymbol{C}(\lambda_i\boldsymbol{E} - \boldsymbol{A})^{-1}\boldsymbol{B}]r_i = w_i\; \forall\; i\; \in\; \{1,2,..., \rho\},
\end{equation}
and left tangential constraints:
\begin{equation}
    l_jH(\mu_j) = l_j[\boldsymbol{C}(\mu_j\boldsymbol{E} - \boldsymbol{A})^{-1}\boldsymbol{B}] = v_j\; \forall\; j\; \in\; \{1,2,..., v\},
\end{equation}
exactly\cite{mayo2007framework}. Where \(\rho + v = N\), or, \(\rho,v = 1, \dots, N\), is considered based on the optimal partition of each data, for further details we refer to see \cite{palitta2022efficient}.
Then the Loewner matrix can be defined as, 
\begin{equation}
    \mathbb{L} = \frac{v_vr_\rho-l_vw_\rho}{\mu_v-\lambda_\rho},
\end{equation}
or, 
\begin{equation}
    \mathbb{L} = \begin{bmatrix} 
                \frac{v_1r_1-l_1w_1}{\mu_1-\lambda_1} & \dots & \frac{v_1r_\rho-l_1w_\rho}{\mu_1-\lambda_\rho} \\
                \vdots & \ddots & \vdots \\
                \frac{v_vr_1-l_vw_1}{\mu_v-\lambda_1} & \dots & \frac{v_vr_\rho-l_vw_\rho}{\mu_v-\lambda_\rho} 
                \end{bmatrix}
                \in \mathbb{C}^{v\times \rho},
\end{equation}

\noindent And the shifted Loewner matrix: 
\begin{equation}
    \mathbb{L}_\sigma = \frac{\mu_vv_vr_\rho-l_vw_\rho}{\mu_v-\lambda_\rho},
\end{equation}
or,
\begin{equation}
    \mathbb{L}_\sigma = \begin{bmatrix} 
                \frac{\mu_1v_1r_1-\lambda_1l_1w_1}{\mu_1-\lambda_1} & \dots & \frac{\mu_1v_1r_\rho-\lambda_\rho l_1w_\rho}{\mu_1-\lambda_\rho} \\
                \vdots & \ddots & \vdots \\
                \frac{\mu_vv_vr_1-\lambda_1l_vw_1}{\mu_v-\lambda_1} & \dots & \frac{\mu_vv_vr_\rho-\lambda_\rho l_vw_\rho}{\mu_v-\lambda_\rho} 
                \end{bmatrix}
                \in \mathbb{C}^{v\times \rho}.
\end{equation}

\noindent Both of these matrices satisfy the Sylvester equations
\begin{equation}
    \mathbf{M}\mathbb{L}-\mathbb{L}\mathbf{\Lambda} = \mathbf{VR}-\mathbf{LW}, \qquad \mathbf{M}\mathbb{L_\sigma}-\mathbb{L_\sigma}\mathbf{\Lambda} = \mathbf{MVR}-\mathbf{LW\Lambda}.
\end{equation}
Then, singular value decomposition of \(\mathbb{L}_\sigma-x\mathbb{L}\) is implemented as,
\begin{equation}
    \mathbf{[Y,\Sigma,X]} = \texttt{svd} (\mathbb{L}_\sigma-x\mathbb{L}), \quad x\in \{s_k\},
\end{equation}
where \(\mathbf{\Sigma}\) is diagonal, \(\mathbf{Y}\) and \(\mathbf{X}\) contain the left and right singular matrices. Based on the first \(r\) dominant modes, the reduced order model is given by,  

\begin{center}
\(\tilde{\mathbf{E}} = -\mathbf{Y}_{r}\mathbb{L}\mathbf{X}_{r} = -\mathbb{L}_{r},\quad \tilde{\mathbf{A}} = -\mathbf{Y}_{r}\mathbb{L}_{\sigma}\mathbf{X}_{r} = -\mathbb{L}_{\sigma r},\quad \tilde{\mathbf{B}} = \mathbf{Y}_{r}\mathbf{V} = \mathbf{V}_{r},\) \\ 
\end{center}

\begin{center}
\(\tilde{\mathbf{C}} = \mathbf{W}\mathbf{X}_r = \mathbf{W}_r,\quad \tilde{\mathbf{D}} = \boldsymbol{0}\). \\
\end{center}

\section{Numerical Results}
\label{NR}

\begin{figure}
\centering
\includegraphics[scale=0.6]{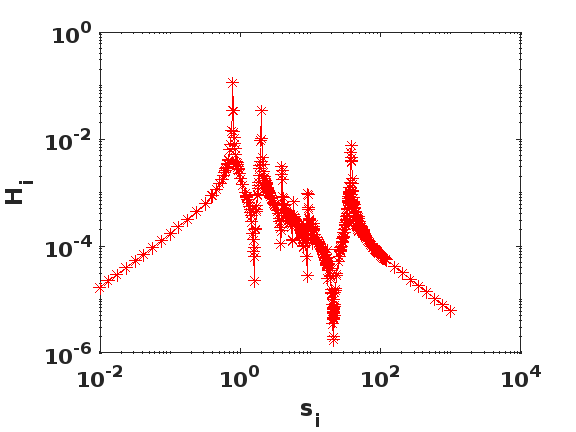} 
\includegraphics[scale=0.6]{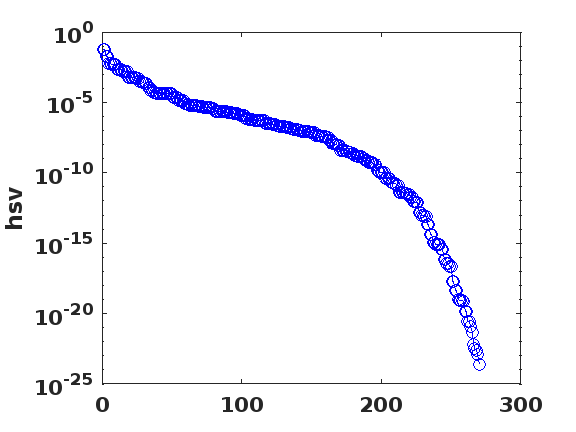}
\caption{System data -- Frequency response from a node (left) and Hankel singular values (right)}%
\label{system_data}
\end{figure}

\begin{figure}
\centering
\includegraphics[scale=0.6]{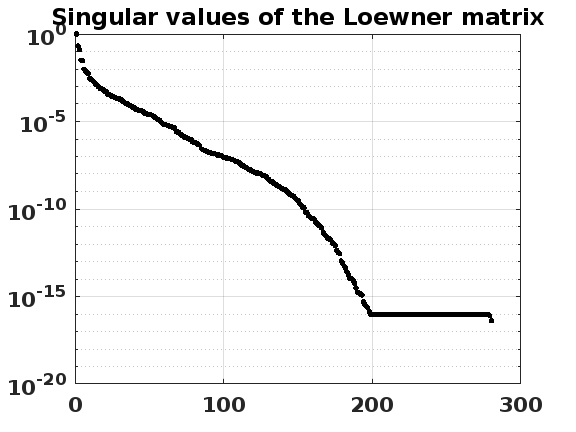} 
\caption{Singular values of the Loewner matrix}%
\label{sing_vals_L}
\end{figure}

\begin{figure}
\centering
\includegraphics[scale=0.5]{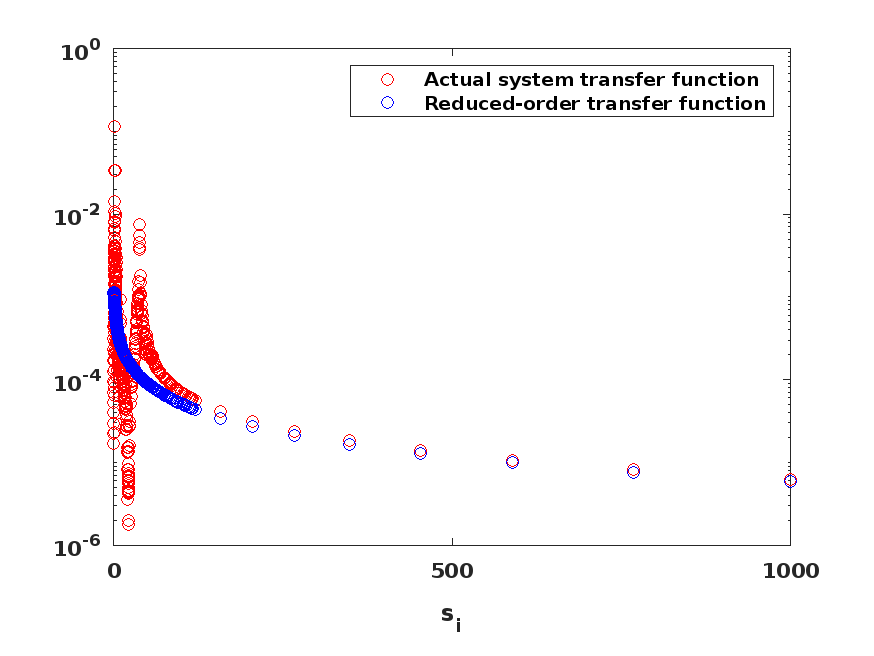}
\includegraphics[scale=0.7]{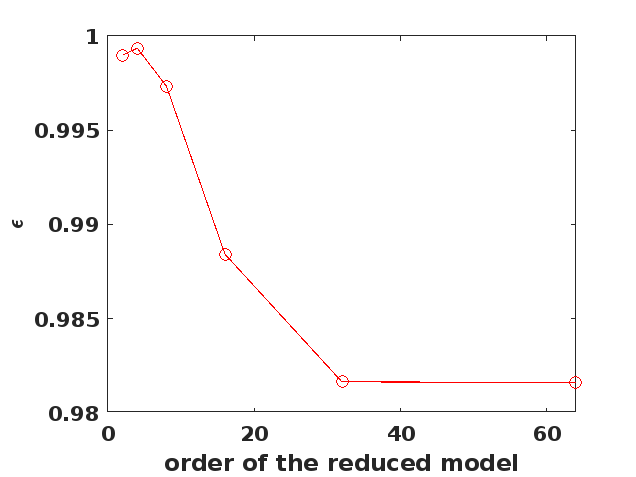}
\caption{Frequency response of the original system and the reduced-order system (\(r = 100\)) (left); Error vs order of the reduced order model (right)}%
\label{compare_performance}
\end{figure}

We demonstrate the Loewner framework on data related to the structural vibrations that occurred as a spaceship was docked at the \textit{International Space Station} \cite{chahlaoui2005benchmark}. The state dimension is \(270\) which implies that the actual system order is 270. This is also evident in the Hankel singular values plot in Figure \ref{system_data} (right). There are 3 inputs and 3 outputs i.e. there are \(9\) nodes in the system. However, we investigate the data from the first node only and it is the transfer function values \(\boldsymbol{H}_i\) at each frequency \(\boldsymbol{s}_{i}\), shown in Figure \ref{system_data} (left). \\ \\
The singular values of the Loewner matrix tell us a lot about where the rank of the reduced order system may lie. We see that the rank can be brought well under 200, Figure \ref{sing_vals_L}. The frequency response data of the original system and the reduced-order model (\(r = 100\)) are plotted in Figure \ref{compare_performance} (left) and it is clear that the reduced-order system is in close agreement with the original one at least in the range of frequencies in consideration. We also observe the error defined as,
\begin{equation}
    \epsilon = \frac{||\boldsymbol{G}_{i} - \boldsymbol{H}_{i}||_{2}}{||\boldsymbol{H}_{i}||_{2}},
\end{equation}
where \(\boldsymbol{G}_{i}\) symbolizes the transfer function values from the reduced order model. Figure \ref{compare_performance} (right) illustrates that the error drops as the order of the reduced system is increased. However, note that it is possible for the error to slightly increase as the order is further increased. 

\newpage

\section{Conclusion and Future Work}
\label{CFW}
In this paper, we employ the Loewner framework on the data obtained from the vibrations in the Russian Service module at the International Space Station. The Loewner model is found to be pretty accurate compared to the original system. About 100 modes are good enough to retain the properties of the actual system. This proves the ability of the Loewner framework to generate a reduced-order model for a given large dynamical system. In the future, we plan to explore how such a framework can be improved to keep the model stable even at higher values of \(r\). 

\newpage

\bibliography{ref}
\bibliographystyle{IEEEtran}

\end{document}